\documentclass[smallextended]{svjour3}       
\smartqed  
\usepackage[ansinew]{inputenc}
\usepackage{amsfonts, color, graphicx}
\usepackage{amssymb,amsmath}
\usepackage{latexsym}
\usepackage{bm}
\usepackage{natbib}
\newtheorem{prop}{Proposition}
\begin{document}
\title{Diffusion Dynamics on the Coexistence Subspace\\ in a Stochastic Evolutionary Game
}

\titlerunning{Diffusive Coexistence in Stochastic Evolutionary Game}        

\author{Lea Popovic         \and
        Liam Peuckert 
}


\institute{L. Popovic  \and
           L. Peuckert \at
              Dept of Mathematics and Statistics, Concordia University, Montreal QC H3G 1M8, Canada\\
              \email{lpopovic@mathstat.concordia.ca}           
}
\date{\today} 
\maketitle

\begin{abstract}
Frequency-dependent selection reflects the interaction between different species as they battle for limited resources in their environment. In a stochastic evolutionary game the species relative fitnesses guides the evolutionary dynamics which fluctuate due to random drift. 
Dependence of species selection advantages on the environment introduces additional possibilities for the evolutionary dynamics. We analyse a simple model in which a random environment allows competing species to coexist for a long time before a fixation of a single species happens. In our analysis we use stability in a linear combination of competing species to approximate the stochastic dynamics of the system by a diffusion on a one dimensional co-existence region. Our method significantly simplifies calculating the probability of first extinction and its expected time, and demonstrates a rigorous model reduction technique for evaluating quasistationary properties of a stochastic evolutionary model.
\end{abstract}

\tableofcontents

\section{Introduction}\label{sec:intro}

Evolutionary games have been useful in modelling reproductive successes of different types of species (individuals) based on their traits and interactions with other types in the population. The fitness of different species types depends on relative proportions (frequencies) of types, on the specific type of interactions, and on the environment. Traditionally these models were deterministic (see \cite{Sigmund98}) and formulated in terms of a system of ordinary differential equations (replicator equation) that describes the evolution of species frequencies. The underlying assumption in deterministic models is that populations are infinite and that the role of fluctuations in determining long term dynamics is negligible. 


Over the last decade stochastic versions of evolutionary games  were used (\cite{Nowak04, Nowak05, Imhof06}, \cite{Traulsen05}) to model the behaviour in finite populations and to reveal the effects caused by fluctuations in the long term dynamics of populations. In the limit as the population size grows some aspects of stochastic models are well approximated by the dynamics of its mean, but the role of fluctuations is crucial for describing events that happen on a very long time scale, such as fixation or extinction. For example, in Moran type models final outcomes are not always determined by the same conditions as in deterministic dynamics. There are examples in which selection completely favours one or the other species for finite values of the population size  where the deterministic dynamics has both species as evolutionarily stable. (For recent reviews of both deterministic and stochastic evolutionary game modelling see \cite{Sandholm07, Hauertrev09}).
 
In stochastic evolutionary models one is particularly interested in chances of reaching different equilibrium or absorption points (extinction and fixation probabilities), as well as the proportion of time needed to reach them (first extinction and fixation times). When the number of types is greater than two both  deterministic and stochastic behaviour becomes much more complex and difficult to quantify. For individual-based (Markov chain) models one needs to solve a system of difference equations which, for large population sizes, is numerically intensive. One can use stochastic differential equations (diffusions) to approximate Markov chain dynamics in large finite populations (e.g.\cite{Traulsen12}) to simplify some of these computations.  For example, calculating the probability of fixation via the optimal stopping theorem would then require solving PDEs with appropriately prescribed boundary values. The validity of a diffusion  approximation depends on the characteristics of the individual-based models as well as on the quantities one wishes to compute. For example, calculating the time to fixation involves events in which the population values are very near the boundary of the state space for the system and the approximation will not be accurate for types whose size is very small. A diffusion approximation is only valid after a careful rigorous justification. 

Specifics of the individual-based model affect the type of diffusion that is appropriate. 
We focus on a model in which the approximating diffusion is degenerate (its fluctuations are negligible in certain directions) and hence the diffusion space is of lower dimension than the original system. A rigorous justification and the derivation of the coefficients for such a diffusion are more complex than in the standard case, but long term behaviour of such degenerate diffusions is much simpler to quantify. Such an  approximation  is a useful tool for model reduction and eliminates the need for time consuming simulations.

There is also an important body of literature examining how random fluctuations due to demographic and environmental stochasticity can contribute to persistence and coexistence in populations during a long transient period of time that precedes the ultimate extinction. (For an excellent survey of recent mathematical results see \cite{Schreiber17}). One can characterize transient meta-stable behaviour of a Markov chain model by quasi-stationary distributions (\cite{Meleardetal09}). When the model is such that stochastic effects can be represented as random perturbation of a dynamical system with a small noise coefficient, then, under additional assumptions on the system, these quasi-stationary distributions concentrate on positive attractors of deterministic dynamics and the probability of extinction decays exponentially with population size (\cite{FaureSchreiber14}). 

The model we consider does not fit into this framework, though it has the same behaviour in its initial time period. However, on a longer time scale our model can with a non-negligible probability (not decaying exponentially with population size) leave the positive attractors of the deterministic dynamics, and can no longer be described as a random perturbation of a dynamical system. 
We complement the aforementioned body of work (\cite{Schreiber17}, \cite{FaureSchreiber14}) by presenting a model whose behaviour in the transient period can be approximated 
by a diffusion of lower dimension (living on a stable manifold of positive attractors for the deterministic dynamics), and by describing a technique that approximates the length of the transient period and the probabilities of extinction. Our model explores the case in which stochasticity is responsible  for long periods of coexistence, but where the proportions of types during coexistence are constrained to a lower dimensional subspace. 

Here we present a simple evolutionary game model in which the environment affects inter-species dynamics. 
The environment is introduced as altering the reproductive fitness of each species type (so the relative fitness of species are different in different  environments). We further let the environment itself be stochastic. Our model considers three species competing in an environment which has two possible states. Two of the species are specialists for the two environments (advantageous in one environment but disadvantaged in the other), and the third is a generalist (indifferent to the environment). We consider a stationary environment (i.i.d. at all reproduction-competition times), and we express long term properties of the model in terms of the environment statistics.
{The effect of the random environment results in long-term stability of a fixed linear combination of the two specialist species, while fluctuating in their proportion relative to the generalist species. This form of stability in a coexistence region allows us to use a degenerate (reduced dimension) diffusion approximating the stochastic dynamics of the system}. This allows us to calculate the first extinction probabilities and expected time to this event, and to quantify the effect of the environment on the long term dynamics.


\section{Stochastic Evolutionary Game Model}

\subsection{Individual based model in a random environment}
Our model is a version of the frequency dependent Moran model with three different species and a constant large population size $N$. Each individual lives for an Exponentially distributed (with rate 1) amount of time at the end of which it is replaced by an individual whose type is chosen at random with chance proportional to its fitness in the current environment. The state of the environment is random from one of two possible states. 

To illustrate the simple model idea let us denote the three species as $C$, $H$ and $M$, and the environment states as {\bf c} and {\bf h}. The two species $C$ and $H$ are specially adapted  to the environments {\bf c} and {\bf h} respectively. The third species $M$ is a generalist which is equally adapted to both of the environments but with less of an advantage than the specialist in either environment. More precisely, $C$ has fitness $1$ in the environment {\bf c} it is specialized for and $0$ in the other environment {\bf h}, $H$ has fitness $1$ in the environment {\bf h} and $0$ in {\bf c}, while $M$ has fitness $1/2$ in both {\bf c} and {\bf h}. The environment has  probability $q$ of being in state {\bf c} and $1-q$ of being in state {\bf h}.

The stochastic dynamics of this evolution game model can be represented by a continuous-time Markov chain $(X_t) _{t\ge 0}=(C_t,H_t,M_t) _{t\ge 0}$ on the state space $\{0,1,\dots,N\}^3$ subject to the constraint of a constant population size $C_t+H_t+M_t=N$, $\forall t\ge 0$. The Markov chain has three absorbing states $(N,0,0), (0,N,0), (0,0,N)$, corresponding to the fixation of any of the species. Outside of the absorbing states the transition kernel is specified as follows. Let \[\mu_c(X)=\frac{C}{N}, \;\mu_H(X)=\frac{H}{N}, \;\mu_M(X)=\frac{M}{N}\] denote the density dependent rates at which individuals of different species are chosen to die, i.e. be replaced, and let \begin{align*}\lambda_C(X)&=q\frac{2C}{(2C+M)}, \;\lambda_H(X)=(1-q)\frac{2H}{(2H+M)}, \\\lambda_M(X)&=q\frac{M}{(2C+M)}+(1-q)\frac{M}{(2H+M)}\end{align*} denote the environment dependent rates at which they are chosen to be born, i.e be the replacement (the fitness takes into account the state of the environment  and the environment specific fitnesst in each state); the possible jumps of $X$ and their rates are given by
\begin{center}
\begin{tabular}{c|c}
$\Delta (C,H,M)$& \mbox{jump rate}\\
\hline
(1,0,-1)&$\lambda_C(X)\mu_M(X)$\\
(1,-1,0)&$\lambda_C(X)\mu_H(X)$\\
(0,1,-1)&$\lambda_H(X)\mu_M(X)$\\
(-1,1,0)&$\lambda_H(X)\mu_C(X)$\\
(0,-1,1)&$\lambda_M(X)\mu_H(X)$\\
(-1,0,1)& $\lambda_M(X)\mu_C(X)$\\
\hline
\end{tabular}
\end{center}
\

\noindent with no change in the process $\Delta (C,H,M)=(0,0,0)$ occurring with rate $\lambda_C(X)\mu_C(X)+\lambda_H(X)\mu_H(X)+\lambda_M(X)\mu_M(X)$.

Due to the constraint  on constant population size the process $X$ is only two dimensional. If $q=\frac12$ there is an obvious symmetry between the two specialists $C$ and $H$. For any $q\in(0,1)$ it is useful to let $D=C-H$ and work with the process $Y=(D,M)$ on $\{-N,\dots, N\}\times\{0,1,\dots,N\}$ subject to the constraints $M+D\le N, M-D\le N$. The process $Y$ is also Markov with absorbing states $(N,0),(-N,0),(0,N)$ and transitions rates expressed explicitely in terms of $D$ and $M$ by 
\begin{center}
\begin{tabular}{c|c}
$\Delta (D,M)$& \mbox{jump rate}\\
\hline
(1,-1)&$\frac{qM(N+D-M)}{N(N+D)}$\\
(2,0)&$\frac{q(N+D-M)(N-D-M)}{2N(N+D)}$\\
(-1,-1)& $\frac{(1-q)M(N-D-M)}{N(N-D)}$\\
(-2,0)&$\frac{(1-q)(N+D-M)(N-D-M)}{2N(N-D)}$\\
(1,1)&$\frac{M(N-D-M)(N+D-2qD)}{2N(N+D)(N-D)}$\\
(-1,1)&$\frac{M(N+D-M)(N+D-2qD)}{2N(N+D)(N-D)}$\\
\hline
\end{tabular}
\end{center}
\

\noindent and no change $\Delta (D,M)=(0,0)$ occurring at rate \[\frac{q[(N+D-M)^2+2M^2]}{2N(N+D)}+\frac{(1-q)[(N-D-M)^2+2M^2]}{2N(N-D)}.\]
The change in the mean is given by
\begin{eqnarray}\label{eq:EDEM}
\begin{aligned}
E[\Delta D_t|Y_t]=&\frac{-(N^2-NM-D^2)(N+D-2qN)}{N(N+D)(N-D)}\\
E[\Delta M_t|Y_t]=&\frac{DM(N+D-2qN)}{N(N+D)(N-D)}\\
\end{aligned}
\end{eqnarray}
and the change in the variance and covariance of the process can be calculated from  

\begin{eqnarray}\label{eq:VarDM}
\begin{aligned}E[(\Delta D_t)^2|Y_t]=&\frac{q(MD+2N^2-2D^2-2MN)}{N(N+D)}\\+&\frac{(1-q)(-MD+N^2-2D^2-2MN)}{N(N-D)}\\
E[(\Delta M_t)^2|Y_t]=&\frac{qM(2N+D-2M)}{N(N+D)}+\frac{(1-q)M(2N-D-2M)}{N(N-D)}\\
E[(\Delta D_t)(\Delta M_t)|Y_t]=&\frac{-qM(-N-2D+M)}{N(N+D)}+\frac{(1-q)M(2N-D-2M)}{N(N-D)}\\
\end{aligned}
\end{eqnarray}

In case $q=\frac12$ there is symmetry in the model between species $C$ and $H$ and their long term coexistence in equal amounts is then expected. Note that  $E[\Delta D_t|Y_t]=0$ when $D=C-H=0$ and that under the same conditions $E[\Delta M_t|Y_t]=0$. This implies that in the symmetric environment case the whole line $\{D=0\}$ is in the null space of the mean dynamics, so no overall push towards a change. In case $q>\frac12$ there is an expected level of dominance of $C$ over $H$. Specifically, when $D=C-H=N(2q-1)$ there is no expected change for $E[\Delta D_t-N(2q-1)|Y_t]=0$ and $E[\Delta M_t|Y_t]=0$. The long term coexistence at this level of difference in their amounts is then expected, as the whole line $\{D=N(2q-1)\}$ is an invariant for the mean. In case $q<\frac12$ there is an expected level of dominance of $H$ over $C$ and the same argument with the roles of $C$ and $H$ reversed show that $\{D=N(1-2q)\}$ is an invariant line for the mean dynamics. 

Importantly, note that on the above invariant lines the variance is not zero, hence the stochastic evolutionary process still changes over time. Before we proceed with a detailed analysis, we first make a few comparisons. 

Consider the version of the model when the environment is not random. If the environment were always in state {\bf h}, i.e. $q=0$, species $C$ would be non-viable and replacement by $H$ would be favoured over replacement by $M$ (by a factor of $2$), hence the dynamics would result in fixation at $H=N, C=M=0$. If the environment were always in state {\bf c}, i.e. $q=1$, then $H$ would be non-viable, replacement by $C$ would be favoured to that by $M$, so the dynamics would result in  fixation at $C=N, H=M=0$. 

One could also consider a corresponding deterministic evolutionary game:  in the environment {\bf c} it would have pay-offs equal to $(2,0,1)$ for $(C,H,M)$ respectively, and  in the environment {\bf h} it would have pay-offs $(0,2,1)$ for $(C,H,M)$ respectively. The replicator equations in a randomized environment with chances $(q,1-q)$ for {\bf (c,h)} predict stable coexistence only when $q=\frac12$. In case $q>\frac12$ the dynamics ultimately always leads to fixation in $C$ and extinction of $H$ and $M$, and in case $q<\frac12$ it leads to fixation in $H$ and extinction of $C$ and $M$.

{Thus, both randomness of the environment and from demographic fluctuations} are required for the dynamics to feature long term coexistence of competing species. In order to quantify its stochastic properties, i.e. evaluate the order of magnitude of its coexistence phase and determine the relative chances of ultimate fixation of competing species, we will derive a diffusion approximation of the individual based model. The  approximation accurately reflects the long term coexistence of the original model, as well as allows us to efficiently approximate the forementioned stochastic quantities.

\subsection{Short time competitive dynamics}

We assume the population size $N$ is large and consider the dynamics of a rescaled process of  $(d,m)=(\frac DN,\frac MN)$ proportions of species in the population, under a factor of $N$ compression of the time scale:  \[y^N_t=(d_t,m_t):=\frac{Y_{Nt}}N=\Big(\frac {D_{Nt}}N,\frac {M_{Nt}}N\Big)\] 
The rescaled process is contained in the triangle \[S=\{-1\le d\le 1, 0\le m\le 1, m+d\le 1, m-d \le 1\}.\]

\begin{figure}
\centering
\includegraphics[totalheight=0.3\textheight,width=1\textwidth]{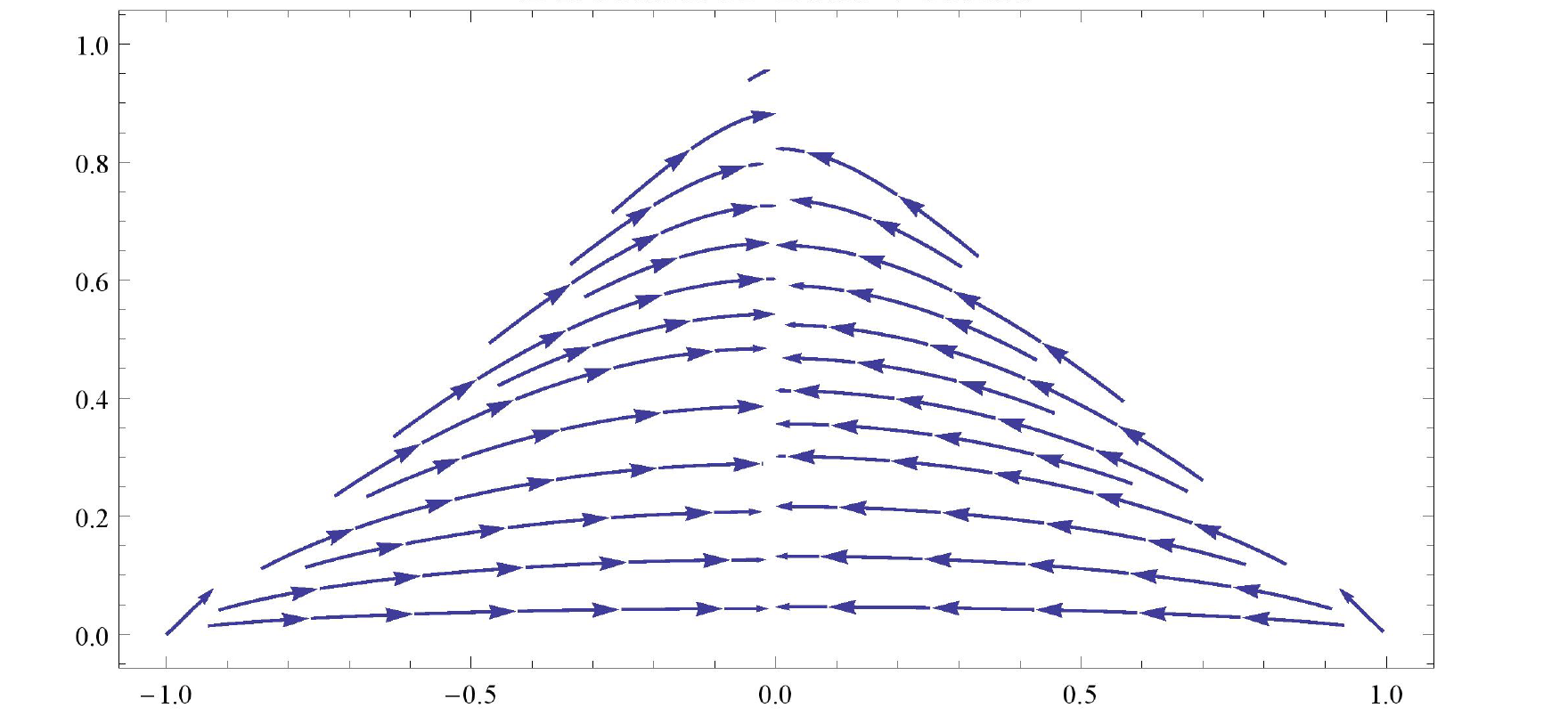}
\caption{\small Gradient of mean change when $q=\frac 12$ (flow of $(b_d,b_m)=E(\Delta d,\Delta m)$ with $d$ on the $x$-axis and $m$ on the $y$-axis).  The flow is directed towards the line $\{d=0\}$ where its gradient vanishes. Also $m$ is increasing (except at $m=0$).}
\label{fig:means}
\end{figure}

\begin{figure}
\centering
\includegraphics[totalheight=0.325\textheight,width=1\textwidth]{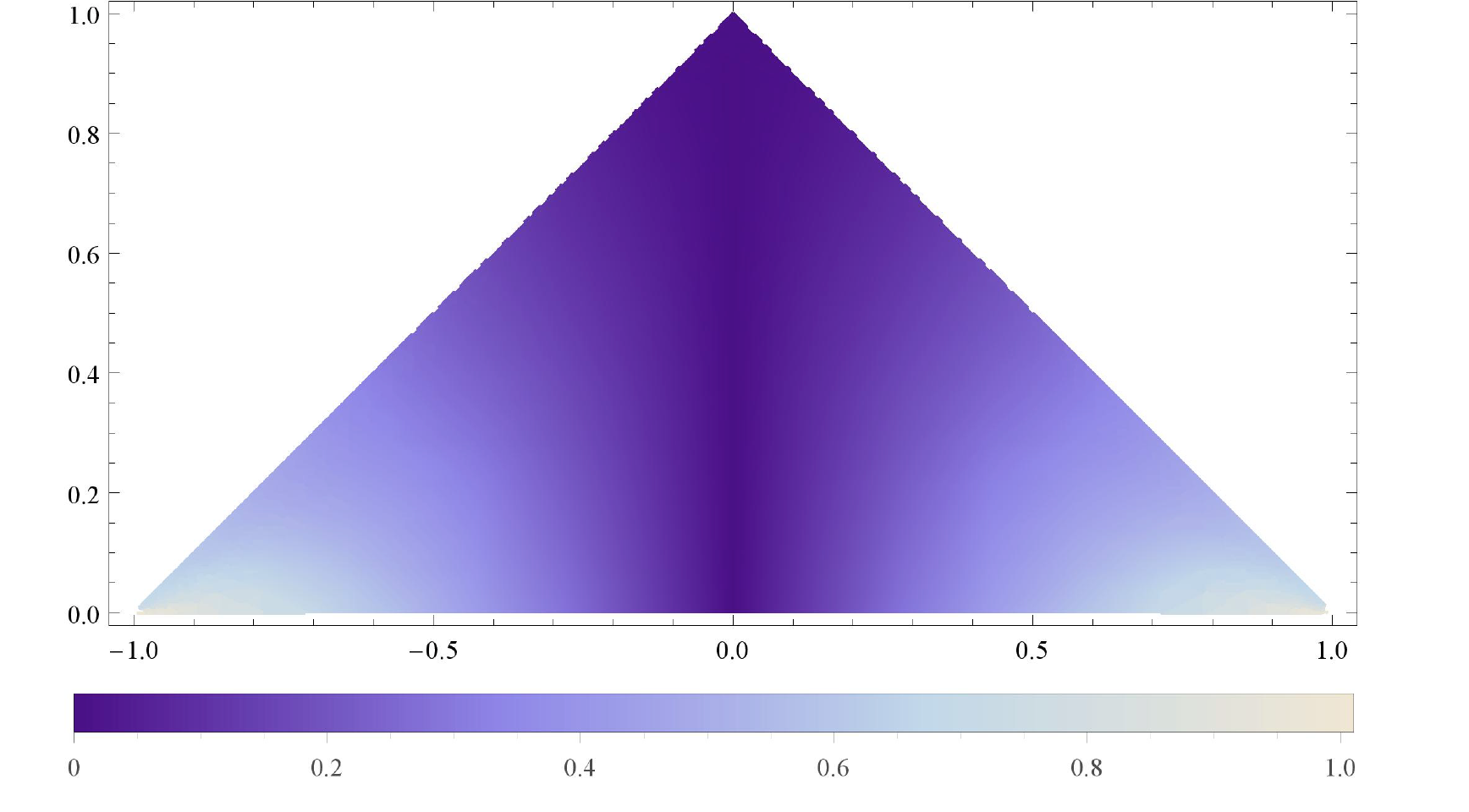}
\caption{Magnitude of the gradient of mean change when $q=\frac 12$ (density plot of $\|(b_d,b_m)\|$ with $d$ on the $x$-axis and $m$ on the $y$-axis). The gradient has the greatest pull at the absorbing points and decreases as it approaches $\{d=0\}$.}
\label{fig:magnitude}\end{figure}

Let the vector $b=(b_d, b_m)$ and matrix $a=((a_{dd}, a_{dm}),(a_{md},a_{mm}))^t$ with $a_{md}\equiv  a_{dm}$ be the functions of $y$ which reflect the infinitesimal change in the mean and (co)variance of the rescaled process:
\[ b_d(y^N_t)=E[\Delta d_t|y^N_t],  b_m(y^N_t)=E[\Delta m_t|y^N_t] \]
\[\frac1Na_{dd}(y^N_t) = E[(\Delta d_t)^2|y^N_t], \frac1Na_{mm}(y^N_t)=E[(\Delta m_t)^2|y^N_t]\]
\[\frac1Na_{dm}(y^N_t)=E[(\Delta d_t)(\Delta m_t)|y^N_t]\]
then, from \eqref{eq:EDEM} and \eqref{eq:VarDM} we have
\begin{eqnarray}\label{eq:moments}
\begin{aligned} 
b_d(y^N_t)=&\frac{-(1-m-d^2)(1+d-2q)}{(1+d)(1-d)}\\
b_m(y^N_t)=&\frac{dm(1+d-2q)}{(1+d)(1-d)}\\
a_{dd}(y^N_t)=&\frac{q(md+2-2d^2-2m)}{(1+d)}+\frac{(1-q)(-md+2-2d^2-2m)}{(1-d)}\\
a_{mm}(y^N_t)=&\frac{qm(2+d-2m)}{(1+d)}+\frac{(1-q)m(2-d-2m)}{(1-d)}\\
a_{dm}(y^N_t)=&\frac{qm(-1-2d+m)}{(1+d)}+\frac{(1-q)m(1-2d-m)}{(1-d)}
\end{aligned}
\end{eqnarray}
The flow of means determined by gradient $(b_d,b_m)$ is shown in Figure~\ref{fig:means} and its magnitude of the mean change is shown in Figure~\ref{fig:magnitude}. 

The scaling of the mean and covariance in the dynamics of the model suggests that on most of the space the mean dominates the  fluctuations. This dominance  $\|b(y)\|\gg \frac{1}{N}\|a(y)\|$ holds only when 
\begin{equation}\label{eq:approxbounds}
|d-(2q-1)|\gg o(1/N)\;\mbox{  and }\;|1-d^2|+|m|\gg o(1/N),
\end{equation} 
while when $|d-(2q-1)|\sim O(1/N)$ or when $\|(d,m)-(\pm 1, 0)\| \sim O(1/N)$ the mean term itself is small and compares in magnitude with the noise. When the mean is dominant the rescaled process $(y^N_t)_{t\ge 0}$ can be approximated by the dynamics of the deterministic solution of the ODE 
\begin{equation}\label{eq:2dode}\bar y_t=\bar y_0+\int_0^tb(\bar y_s)ds
\end{equation}
The following result insures that away from the three absorbing points the deterministic approximation holds for $y^N_t$ on any finite time interval. 

\begin{prop}\label{prop:ODE} If \mbox{$y^N_t=\frac{Y_{Nt}}N$} is the rescaled process started from $y^N_0$ and $\bar y$ is the deterministic solution of \eqref{eq:2dode} started from \mbox{$\bar y_0=(d_0,m_0)$}, with \mbox{$y^N_0\Rightarrow \bar y_0$}, then for any $0<\varepsilon<1$, $N\ge 1$, $T>0$ there is a constant  $\sup_NL(N,T)<\infty$, such that for  $1-d_0^2>\varepsilon$, $0<m_0<1$, and $\tau^N_\epsilon=\inf\{t:1-(d^N_t)^2 \le\epsilon\}$, we have
\begin{equation*}\label{eq:propodeapprox}
E(\sup_{t\le T\land \tau^N_\epsilon} \|y^N_t-\bar y_t\|^2)\le \frac{L(N,T)^2}{{N}}.
\end{equation*}
\end{prop}
The proof of this result follows from the approximation of density dependent Markov processes by ODE derived by \cite{Kurtz78}, and relies on the bounded behaviour of the rescaled jumps and the Lipschitz continuity of the vector field of the limiting ODE) (For more details on theory of approximating density dependent processes see \cite{Kurtz87}). 

The evolution of the deterministic dynamics given by \eqref{eq:2dode} on the state space  $S=\{-1\le d\le 1, 0\le m\le 1, m+d\le 1, m-d \le 1\}$ is as follows. 
The gradient vector $b=(b_d,b_m)$ has a component $b_d$ vanishing iff either $d=2q-1$ or $(d, m)\in\{(1, 0),(-1,0)\}$ and the other component $b_m$ vanishing iff either $d\in\{0,2q-1\}$ or $m=0$.
The line $\Gamma=\{(2q-1,m)\}$ is stable, and since the function $\varphi(y)=(d-2q+1)^2$ is a Lyapunov function for this dynamical system \[\nabla \varphi\cdot b=-2(1-m-d^2)(1+d-2q)^2/(1-d)(1+d)\le 0,\]
$\Gamma$ is also globally attracting on $S$ excluding the corners of this triangle $(d, m)\in\{(1, 0),(-1,0),(0,1)\}$. These corner points are also fixed points for the dynamical system \eqref{eq:2dode} (and absorption points for the original stochastic process).
 In addition $\nabla m= b_m\ge 0$ on the whole triangle, so  $m$ is always increasing (except on the line $\{m=0\}$). The magnitude of $b$ is highest in the neighbourhood of the two corners $(d, m)\in\{(1, 0),(-1,0)\}$ and lowest near the line $d=2q-1$, see Figures~\ref{fig:means} and~\ref{fig:magnitude}.
 

The trajectory from any initial point in $S$ except the forementioned corners evolves by converging to the stable line $\Gamma$. Let $\Phi(d,m)$ denote the map that takes a point $y=(d,m)$ to $\Gamma$, that is, the intersection of the trajectory of $\bar y$ started from $\bar y_0=(d,m)$ with the line $\Gamma$. Let $\boldsymbol{m}^*(d,m)$ be the $m$-coordinate of the map, so $\Phi(d,m)=(2q-1,\boldsymbol{m}^*)$.
From \eqref{eq:2dode} and \eqref{eq:moments} we have an autonomous system 
\begin{eqnarray*}\label{eq:curveode}
\frac{\partial d}{\partial m}=\frac{-(1-m-d^2)}{dm},\quad  \frac{\partial m}{\partial d}=\frac{dm}{-(1-m-d^2)}
\end{eqnarray*}
solved by
\begin{eqnarray*}\label{eq:curveeq}
d(m)&=\pm \sqrt{1-2m+Cm^2},\quad m(d)=\frac{1-d^2}{1+\sqrt{1-C(1-d^2)}}
\end{eqnarray*}
where the constant $C$ is determined by initial conditions $y_0=(d,m)$ to be 
\[C=\frac{-1+2m+d^2}{m^2}.\]
Since $\Gamma=\{d=2q-1\}$ we get that the $m$-coordinate of the point of intersection of the trajectory starting from $\bar y_0=(d,m)$ with $\Gamma$ is 
\begin{eqnarray}\label{eq:mstar}
\boldsymbol{m}^*
=\frac{4q(1-q)}{1+\sqrt{1-\frac{-1+2m+d^2}{m^2}4q(1-q)}}
\end{eqnarray}

Proposition~\ref{prop:ODE} implies that the first phase of the model dynamics is governed by a strong competitive drive between the two specialists $C$ and $H$ towards their coexistence proportions $c-h=2q-1$ determined by the distribution of the two environmental states. The proportion of the generalist is $\boldsymbol m^*$ when the coexistence domain is reached and depends on the initial proportions of the three species. The rest of the time  stochastic dynamics will continue in a close neighbourhood of the domain of coexistence $\Gamma$, where, as we next prove, fluctuations in the proportion of the generalist versus the two specialists will persist for a very long time. 

\subsection{Long time coexistence and persistent fluctuations}

In the initial time period the species count is well approximated by the associated deterministic curve. Asymptotic stability of $\Gamma$ for this dynamics would suggest that in the long run the species count would converge to $\Phi(d,m)$ on $\Gamma$. However, \eqref{eq:approxbounds}  shows that close to $\Gamma$ this approximation is no longer accurate, hence we need to examine the original process on a longer time scale .
Since in a neighbourhood of the stable line $\Gamma$ fluctuations are non-negligable, a diffusion approximation is necessary. 

The rescaled process $y_t^{N^2}=y^N_{Nt}=Y_{N^2t}/N$ under a factor of $N^2$ compression of the time scale behaves like the two-dimensional diffusion  
\begin{equation}\label{eq:2ddiff}
\tilde y_t=\tilde y_0+N \int_0^t b(\tilde y_s) ds + \int_0^t\sigma (\tilde y_s) dW_s
\end{equation}
where $W=(w^d,w^m)$ is a standard Brownian motion in $\mathbb R^2$, $b(y)$ is as in \eqref{eq:moments} and $\sigma(y)$ is the matrix such that $\sigma^t(y)\sigma(y)=a(y)$ from \eqref{eq:moments}.
This diffusion has a very strong mean and fluctuations only of size $O(1)$. When $N$ is large and the process is at a point outside of $\Gamma$, the very large mean term overpowers the noise and quickly carries the process onto the stable line $\Gamma$. Once the process is on $\Gamma$ the mean term is zero so the noise can move the process away, but as soon as its distance from $\Gamma$ is more than $O(1/N)$ the mean carries it back, this time to another point on $\Gamma$. The noise acts as a perpetual source of perturbations for the stabilizing deterministic dynamics. In the limit of a large population size this becomes equivalent to a diffusive process on $\Gamma$. 

In order to calculate the characteristics of the limiting process on $\Gamma$ we need to use the map $\Phi$ which projects each point to its trajectory destination on $\Gamma$ after it is perturbed by the diffusion \eqref{eq:2ddiff}. Recall that $\Phi(d,m)=(2q-1,\boldsymbol{m}^*)$ where $\boldsymbol{m}^*$ is a function of $(d,m)$ through the constant $C=C(d,m)$ (see  \eqref{eq:mstar} above)
\[\boldsymbol{m}^*=\frac{4q(1-q)}{1+\sqrt{1-C4q(1-q)}}.\] 
Considering $\Phi(d,m)=(2q-1,\frac{4q(1-q)}{1+\sqrt{1-C(d,m)4q(1-q)}})$ as a function applied to the two-dimensional diffusion \eqref{eq:2ddiff} we can use Ito's lemma in order to determine the stochastic behaviour of $\Phi(d_t,m_t)$. For this we will need the following derivative calculations:
\begin{align}\label{eq:partialmstar}
&\frac{\partial \boldsymbol{m}^*}{\partial d}=\frac{(\boldsymbol{m}^*)^3}{2(4q(1-q)-\boldsymbol{m}^*)}\frac{\partial C}{\partial d}\\
&\frac{\partial \boldsymbol{m}^*}{\partial m}=\frac{(\boldsymbol{m}^*)^3}{2(4q(1-q)-\boldsymbol{m}^*)}\frac{\partial C}{\partial m}\nonumber\\
&\frac{\partial^2 \boldsymbol{m}^*}{\partial^2 d}=\frac{(\boldsymbol{m}^*)^3}{2(4q(1-q)-\boldsymbol{m}^*)}\frac{\partial^2 C}{\partial^2 d}+\frac{(\boldsymbol{m}^*)^2(12q(1-q)-2\boldsymbol{m}^*))}{2(4q(1-q)-\boldsymbol{m}^*)^2}\frac{\partial \boldsymbol{m}^*}{\partial d}\frac{\partial C}{\partial d}\nonumber\\
&\quad\quad\;\;=\frac{(\boldsymbol{m}^*)^3}{2(4q(1-q)-\boldsymbol{m}^*)}\frac{\partial^2 C}{\partial^2 d}+\frac{(\boldsymbol{m}^*)^5(12q(1-q)-2\boldsymbol{m}^*))}{4(4q(1-q)-m^*)^3}\Big(\frac{\partial C}{\partial d}\Big)^2\nonumber\\
&\frac{\partial^2 \boldsymbol{m}^*}{\partial^2 m}=\frac{(\boldsymbol{m}^*)^3}{2(4q(1-q)-\boldsymbol{m}^*)}\frac{\partial^2 C}{\partial^2 m}+\frac{(\boldsymbol{m}^*)^2(12q(1-q)-2\boldsymbol{m}^*))}{2(4q(1-q)-\boldsymbol{m}^*)^2}\frac{\partial \boldsymbol{m}^*}{\partial m}\frac{\partial C}{\partial m}\nonumber\\
&\quad\quad\;\;=\frac{(\boldsymbol{m}^*)^3}{2(4q(1-q)-\boldsymbol{m}^*)}\frac{\partial^2 C}{\partial^2 d}+\frac{(\boldsymbol{m}^*)^5(12q(1-q)-2\boldsymbol{m}^*))}{4(4q(1-q)-\boldsymbol{m}^*)^3}\Big(\frac{\partial C}{\partial m}\Big)^2\nonumber\\
&\frac{\partial^2 \boldsymbol{m}^*}{\partial d\partial m}=\frac{(\boldsymbol{m}^*)^3}{2(4q(1-q)-\boldsymbol{m}^*)}\frac{\partial^2 C}{\partial d\partial m}+\frac{(\boldsymbol{m}^*)^2(12q(1-q)-2\boldsymbol{m}^*))}{2(4q(1-q)-\boldsymbol{m}^*)^2}\frac{\partial \boldsymbol{m}^*}{\partial m}\frac{\partial C}{\partial d}\nonumber\\
&\quad\quad\;\;=\frac{(\boldsymbol{m}^*)^3}{2(4q(1-q)-\boldsymbol{m}^*)}\frac{\partial^2 C}{\partial d\partial m}+\frac{(\boldsymbol{m}^*)^5(12q(1-q)-2\boldsymbol{m}^*))}{4(4q(1-q)-\boldsymbol{m}^*)^3}\Big(\frac{\partial C}{\partial d}\frac{\partial C}{\partial m}\Big)\nonumber
\end{align}
and since $C(d,m)=(-1+2m+d^2)/m^2$ we also have that
\begin{eqnarray}\begin{aligned}\label{eq:partialC}
&\frac{\partial C}{\partial d}=\frac{2d}{m^2},\\
&\frac{\partial C}{\partial m}=\frac{2(1-d^2-m)}{m^3}\\
&\frac{\partial^2 C}{\partial^2 d}=\frac{2}{m^2},\\
&\frac{\partial^2 C}{\partial^2 m}=-\frac{2(3-3d^2-2m)}{m^4},\\
&\frac{\partial^2 C}{\partial d\partial m}=-\frac{4d}{m^3}
\end{aligned}
\end{eqnarray}

\

Ito's formula applied to $\Phi(d,m)=(2q-1,\boldsymbol{m}^*(d,m))$ implies that  the second coordinate of $\Phi$ is a one-dimensional diffusion process. Its mean coefficient is given by
\begin{align}\label{eq:diffonGbeta}
\beta_{\boldsymbol{m}^*}(d,m)&=N\Big(\frac{\partial \boldsymbol{m}^*}{\partial d}b_d+\frac{\partial \boldsymbol{m}^*}{\partial m}b_m\Big)+\frac12\Big(\frac{\partial^2 \boldsymbol{m}^*}{\partial d^2}a_{dd}+\frac{\partial^2 \boldsymbol{m}^*}{\partial m^2}a_{mm}+2\frac{\partial^2 \boldsymbol{m}^*}{\partial d\partial m}a_{dm}\Big)\nonumber\\
&=\frac12\frac{\partial^2 \boldsymbol{m}^*}{\partial d^2}a_{dd}+\frac12\frac{\partial^2 \boldsymbol{m}^*}{\partial m^2}a_{mm}+\frac{\partial^2 \boldsymbol{m}^*}{\partial d\partial m}a_{dm}
\end{align}
since by \eqref{eq:partialmstar}, \eqref{eq:partialC} and \eqref{eq:moments} we have that
\begin{align*}
&\frac{\partial \boldsymbol{m}^*}{\partial d}b_d+\frac{\partial \boldsymbol{m}^*}{\partial m}b_m\\
&=\frac{(\boldsymbol{m}^*)^3}{2(4q(1-q)-\boldsymbol{m}^*)}\left(\frac{2d}{m^2}\frac{-(1-m-d^2)(1+d-2q)}{(1+d)(1-d)}+\frac{2(1-m-d^2)}{m^2}\frac{dm(1+d-2q)}{(1+d)(1-d)}\right)\\&=0,
\end{align*}
Its diffusion coefficient is calculated by combining the contributions to fluctuations from the two independent noise components in \eqref{eq:2ddiff} given by
\begin{align}\label{eq:diffonGalpha}
\alpha_{\boldsymbol{m}^*}(d,m)&=\sqrt{\Big(\frac{\partial \boldsymbol{m}^*}{\partial d}\sigma_{dd}+\frac{\partial \boldsymbol{m}^*}{\partial m}\sigma_{md}\Big)^2+\Big(\frac{\partial \boldsymbol{m}^*}{\partial d}\sigma_{dm}+\frac{\partial \boldsymbol{m}^*}{\partial m}\sigma_{mm}\Big)^2}\nonumber\\
&=\sqrt{\Big(\frac{\partial \boldsymbol{m}^*}{\partial d}\Big)^2a_{dd}+\Big(\frac{\partial \boldsymbol{m}^*}{\partial m}\Big)^2a_{dd}+2\Big(\frac{\partial \boldsymbol{m}^*}{\partial d}\Big)\Big(\frac{\partial \boldsymbol{m}^*}{\partial m}\Big)a_{dm}}
\end{align}
and consequently 
\begin{equation}\label{eq:1ddiffonG}
\boldsymbol{m}^*_t=\boldsymbol{m}^*_0+\int_0^t \beta_{\boldsymbol{m}^*}(d,m)ds + \int_0^t\alpha_{\boldsymbol{m}^*}(d,m)dw_s, \mbox{ for }\boldsymbol{m}^*\in [0,1]
\end{equation}
where $w$ is a standard Brownian motion in $\mathbb{R}^1$.

The large component in the mean disappears confirming our earlier claim that the limit is a regular diffusion  whose state space is the stable line $\Gamma$. This is due to the fact that the contribution to the mean in the perpendicular direction to $\Gamma$ of $\nabla \Phi\cdot b=0$ vanishes, as there is no change in the projection map $\Phi$ along the flow of the mean. The only component in the mean that remains is from the change in $\Phi$ in the direction tangential to $\Gamma$.
A rigorous statement for this one-dimensional approximation of the long term behaviour of our original stochastic process is given by the following result.

\begin{prop}\label{prop:SDE1} If \mbox{$y^{N^2}_t=\frac{Y_{N^2t}}{N}$} is the rescaled process started at a point $y^{N^2}_0$ in a neighbourhood $\|y^{N^2}_0-\Gamma\|\le N^{-\delta}, \delta\in(0,\frac12)$ of the stable line, 
and  if \mbox{$\tau=\inf\{t:m_t\in\{0,1\}\}$} is the time to extinction of the first species in the process; then as $N\to \infty$ the rescaled process stays in the neighbourhood  with 
$P(\sup _{t\le \tau}\|y^{N^2}_t-\Gamma\|\leq N^{-\delta}) \to 1,$ 
and we have convergence of the stopped process: \[y^{N^2}_{\cdot\land \tau}\Longrightarrow \boldsymbol{m}^*_{\cdot\land \tau},\] to the diffusion on $\Gamma$ given by \eqref{eq:1ddiffonG} stopped when it first reaches its boundary point $\{0,1\}$.
\end{prop}

The proof of this result is a consequence of theorem for convergence to degenerate diffusions derived by \cite{Katz91}. (A similar technique was used in \cite{DurPop09} to establish approximation by a degenerate diffusion for a model of genetic subfunctionalization, and an excellent practical exposition is given in \cite{Plot17}). It relies on the strong stability of the line $\Gamma$, differentiability of the projection map $\Phi$ onto $\Gamma$, as well as on the good behaviour of the jumps of the rescaled Markov process. One first uses Lyapunov function $\varphi$  to show that the distance of the process $y^{N^2}$ to $\Gamma$ converges to $0$, followed by a convergence result for integrals with respect to semi-martingales in order to identify that the limiting process on $\Gamma$ is indeed given by the Ito's lemma calculation. Note that we could have combined Proposition~\ref{prop:SDE1} with Proposition~\ref{prop:ODE}  to include the initial period for the process $y^{N^2}$ by starting it anywhere in the domain of attraction $S\setminus\{(-1,0),(1,0)\}$ of $\Gamma$ (this is described in the proof in the Appendix).

Proposition~\ref{prop:SDE1} implies that the second, much longer, phase of the model dynamics consists of trade-off between a fluctuating proportion of the generalist $M$ and a proportion of a combination of  specialists $C$ and $H$. During this time the combination of the latter two is such that their difference is kept approximately at its 'coexistence balance' (as dictated by the distribution of the random environment). Given that this period lasts for a time of order $N^2$ one could try to find its approximate distribution on the line $\Gamma$, using the law of diffusion $\boldsymbol{m}^*_t$ conditioned on $\tau>t$. In this manner one can estimate the occupation time distribution during the long transient period of coexistence. For deriving quasi-staitonary distributions see \cite{Martinez95} Theorem B, \cite{Meleardetal09},Theorem 5.2. 
Ultimately this phase ends with either the fixation of species $M$ or their extinction. In the latter case only $C$ and $H$ remain and continue a stochastic competition until ultimate fixation.

If the process loses all $M$ species first, at the above stopping time $\tau$ the original process $y$ is at  $y_{\tau}=(2q-1,0)$, then the only possible transitions forward in time are $\Delta y\in\{(2,0), (-2,0)\}$ which keep $y$ on the line $\{m=0\}$ while changing the 'coexistence balance' of $c$ and $h$ away from $\{d=2q-1\}$. The probability of fixation (when not at $M$) and the expected time to fixation is then determined by this phase where the process is well approximated by a diffusion $\boldsymbol{d}^*$ with a small noise coefficient. We can make an approximation by $y^{N^2}$ started at $y_0=(d_0,0)$ on the line $\{m=0\}$ in terms of a mean reverting Gaussian process with small noise: $\boldsymbol{d}^*_t=\boldsymbol{d}^*_0+\int_0^t\beta_{\boldsymbol{d}^*}(d^*_s)ds+\frac{1}{\sqrt{N}}\int_0^t\sqrt{\alpha_{\boldsymbol{d}^*}(d^*_s)}dw_s$,  for $\boldsymbol{d}^*\in [-1,1]$ with $\beta_{\boldsymbol{d}^*}(d)=-(d-(2q-1))$ and $\alpha_{\boldsymbol{d}^*}(d)=2-2d(2q-1)$. Due to the limited space constraint, the noise of this process after a long time leads to ultimate fixation. 

\subsection{Probabilities of fixation and its expected time}
To illustrate the usefulness of our results we calculate  the probabilities of fixation and time to fixation from the derived processes.
In the original rescaled process $y^{N^2}=Y_{N^2t}/N\,$ let $\tau_c,\tau_h,\tau_m$ denote the times to extinction of each of the species
\begin{align*}\tau_C=\inf\{t: c_t=0\},\, \tau_H=\inf\{t: h_t=0\},\,\tau_M=\inf\{t: m_t=0\}\end{align*} 
where we assume $\inf\{\emptyset\}=\infty$, let 
$\tau^e$ denote the time of the first extinction and $\tau^f$  the time of fixation, i.e. second extinction. Note, $\tau^e<\infty$ and $\tau^f<\infty$, and it is possible for $\tau^e=\tau^f$. 
Denote the probabilities of first extinction by 
\begin{align*}p^e_C=P(\tau^e=\tau_C),\, p^e_H=P(\tau^e=\tau_H),\, p^e_M=P(\tau^e=\tau_M)\end{align*}
and the probabilities of fixation by
\begin{align*}p^f_C=P(y^{N^2}_{\tau^f}=C),\, p^f_H=P(y^{N^2}_{\tau^f}=H),\, p^f_M=P(y^{N^2}_{\tau^f}=M).\end{align*}

Proposition~\ref{prop:ODE} implies that as long as the starting value is away from the corners $(d, m)\in\{(1, 0),(-1,0), (2q-1,1)\}$ the time to reach the 'coexistence balance' line $(2q-1,m)$ is negligible on the $N^2$ time scale. 

Proposition~\ref{prop:SDE1} further implies that $\tau^e$ is well approximated from the stopping time $\tau=\inf\{t:\boldsymbol{m}^*_t\in\{0,1\}\}$ on the diffusion $\boldsymbol{m}^*$:\\
$\cdot$ if  $\boldsymbol{m}^*_\tau=1$ then both $\tau_C=\tau_H$ and $\tau^e=\tau^f$ are well approximated by $\tau$;\, both  $p^e_{C+H}$ and $p^f_M$ are well approximated by $P(\boldsymbol{m}^*_\tau=1)$;\\
$\cdot$  if  $\boldsymbol{m}^*_\tau=0$ then $\tau^e=\tau_M$ is well approximated by $\tau$;\, both $p^e_M$ and $p^f_C+p^f_H$ are well approximated by $P(\boldsymbol{m}^*_\tau=0)$. \\
So, letting $p_M=P(\boldsymbol{m}^*_\tau=0)$ we have: 
\begin{eqnarray}\label{eq:timesprobs}
\begin{aligned} 
E(\tau^e)&\approx E(\tau),\; E(\tau^f)\approx E(\tau| \mbox{first extinction is not M})\\
p^e_{C+H}&=1-p^e_M\approx p_M,\; p^f_M\approx 1-p_M
\end{aligned}
\end{eqnarray}

Since $\boldsymbol{m}^*$ is a one-dimensional diffusion the above quantities have analytical formulae that can certainly be explicitly computed.
For a  one-dimensional diffusion of the form \[dx_t=\beta(x_t)dt+\sqrt{\alpha(x_t)}dw_t\] we can use the natural scale function (unique up to linear map)
\[\phi(x)=\int^x \exp\{\int^y \frac{-2\beta(z)}{\alpha(z)}dz\}dy\]
from which started from a point $x\in[x_0,x_1]$ the hitting time probabilities $\tau_{x_0}$ and $\tau_{x_1}$ of any two points $x_0<x_1$ are then computed as 
\[P_x(\tau_{x_0}<\tau_{x_1})=\frac{\phi(x_1)-\phi(x)}{\phi(x_1)-\phi(x_0)}.\] 
We can also use the natural scale function to calculate the expected hitting time $E_x[\tau_{x_0}\land \tau_{x_1}]$ started from a point $x\in[x_0,x_1]$ via the Green's function on the interval $[x_0,x_1]$
\begin{align*}G(x,y)&=2\frac{(\phi(x_1)-\phi(x))(\phi(y)-\phi(x_0))}{\phi(x_1)-\phi(x_0)}\frac{1}{\phi'(x)\alpha(x)},y\le x;\\
&=2\frac{(\phi(x)-\phi(x_0))(\phi(x_1)-\phi(y))}{\phi(x_1)-\phi(x_0)}\frac{1}{\phi'(x)\alpha(x)},y> x;\end{align*}
from which the expected hitting time is computed as 
\[E_x[\tau_{x_0}\land \tau_{x_1}]=\int_{x_0}^{x_1} G(x,y)dy.\]

For ease of display we calculate explicitly the formulae above in the symmetric environment case $q=\frac12$. First, for the coefficients of the two-dimensional process \eqref{eq:moments} we get 
\begin{align*}
a_{dd}=\frac{2-2m-2d^2-md^2}{1-d^2}, \; a_{mm}=\frac{m(2-2m-d^2)}{1-d^2},\; a_{md}=\frac{-dm(1+m)}{1-d^2}
\end{align*}
Then for $\boldsymbol{m}^*$ the partial moments \eqref{eq:partialmstar} the mean and variance formulae \eqref{eq:diffonGbeta} and \eqref{eq:diffonGalpha}, 
after some algebraic simplifications, give the following expressions for coefficients:
\begin{align*}
\beta_{\boldsymbol{m}^*}(d,m)&=\frac{-m(-2d^4-2(m-1)(-2+3\sqrt{-d^2+(m-1)^2}+3m))}{2(-1+d^2)\sqrt{-d^2+(m-1)^2}(\sqrt{-d^2+(m-1)^2})^3}\\& +\frac{-m(d^2(6-4m+3m\sqrt{-d^2+(m-1)^2}+3m^2))}{2(-1+d^2)\sqrt{-d^2+(m-1)^2}(\sqrt{-d^2+(m-1)^2})^3}\\
\alpha_{\boldsymbol{m}^*}(d,m)&=-\frac{m(-2+d^2+2m)}{(\sqrt{-d^2+(-1+m)^2}+m)^4}
\end{align*}
Note from above that the diffusion inherently stops when the $m=0$ boundary point is reached, but is not defined when $m=1$. The first case corresponds to the event of loss of the generalist species $M$, while the second case corresponds to fixation in $M$.

The mapping of the starting value $(d_0,m_0)$ onto the line $\Gamma=\{d=0\}$ is given by 
\[\boldsymbol{m}^*(d_0,m_0)=\frac{1}{1+\sqrt{1-\frac{d_0^2+2m_0-1}{m_0^2}}}\]
Using $x=\boldsymbol{m}^*(d_0,m_0)$ the probability $p_M$ is calculated as
\[p_M=P_x(\boldsymbol{m}^*_\tau=0)=\frac{\phi(1)-\phi(x)}{\phi(1)-\phi(0)}\]  
from the natural scale function $\phi(x)=\int^x \exp\{\int^y \frac{-2\beta(z)}{\alpha(z)}dz\}dy$ by numerically integrating 
the inner integral 
over the area projected by the ODE onto the line $\Gamma=\{d=0\}$ between $m=y$ and $m=0$. This corresponds to the area defined by
\[-\sqrt{1-2m+\frac{-1+2y}{y^2}m^2}\le d\le \sqrt{1-2m+\frac{-1+2y}{y^2}m^2}, 0\le m\le y\]
Note that when $y=1$ this corresponds to the integral over the whole triangle $S$.
The expected time to first extinction is then calculated as \[E_x(\tau)=\int_0^1 G(x,y)dy\] by further numerical integration to obtain all values the scale function $\phi(y)$ over $y\in[0,1]$ to be used in the expression
\begin{align*}G(x,y)&=2\frac{(\phi(1)-\phi(x))(\phi(y)-\phi(0))}{\phi(1)-\phi(0)}\frac{1}{\phi'(x)\alpha(x)},y\le x;\\
&=2\frac{(\phi(x)-\phi(0))(\phi(1)-\phi(y))}{\phi(1)-\phi(0)}\frac{1}{\phi'(x)\alpha(x)},y> x;\end{align*}

Explicit values from the one dimensional diffusion for the extinction probability $p_M$ and the expected extinction time $E_{y_0}(\tau)$ are shown alongside simulation results in Figures~\ref{fig:pm} and~\ref{fig:Etau} below. The explicit values are based on numerical integration for calculating $\phi$ at $1000$ points. The graph of $p_M$ and of $E_{y_0}(\tau)$ as functions of the starting point for the process $y_0=(d_0=0,m_0)$ is made by interpolation between computed values for $1000$ different starting points.


\subsection{Simulations results}
\begin{figure}
\centering
\includegraphics[totalheight=0.325\textheight,width=0.75\textwidth]{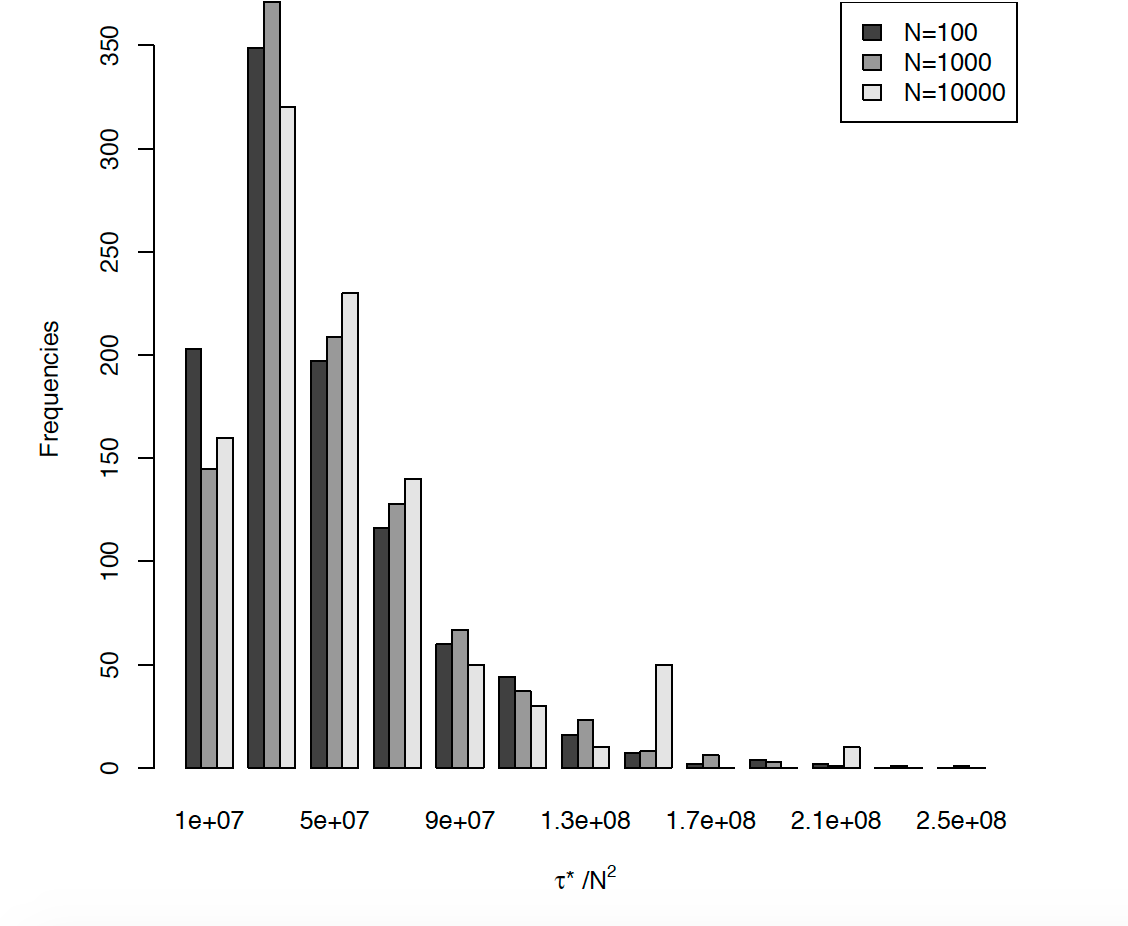}
\caption{\small Distribution of the time to extinction $\tau^e$  for $Y$ starting from $y_0=(0,\frac13)$ for three different values of total population size $N=10^2, 10^3, 10^4$  ($\tau^e$ is scaled by $N^2$ in each case).}
\label{fig:41}
\end{figure}

We compared the results based on the diffusion approximations with simulations of the original stochastic model $Y$ in the symmetric environment $q=\frac12$. 
To assess the sensitivity to the population size we simulated the original process with different orders of magnitude for the population size: $N=100,1000,10000$ and estimated the probability distribution of the time to extinction $\tau^e$ based on outcomes of $1000$ simulation runs for each choice of $N$. The histograms in Figure~\ref{fig:41} show remarkable similarity, and a heavy tail in their distribution 

\begin{figure}
\centering
\includegraphics[totalheight=0.325\textheight,width=0.75\textwidth]{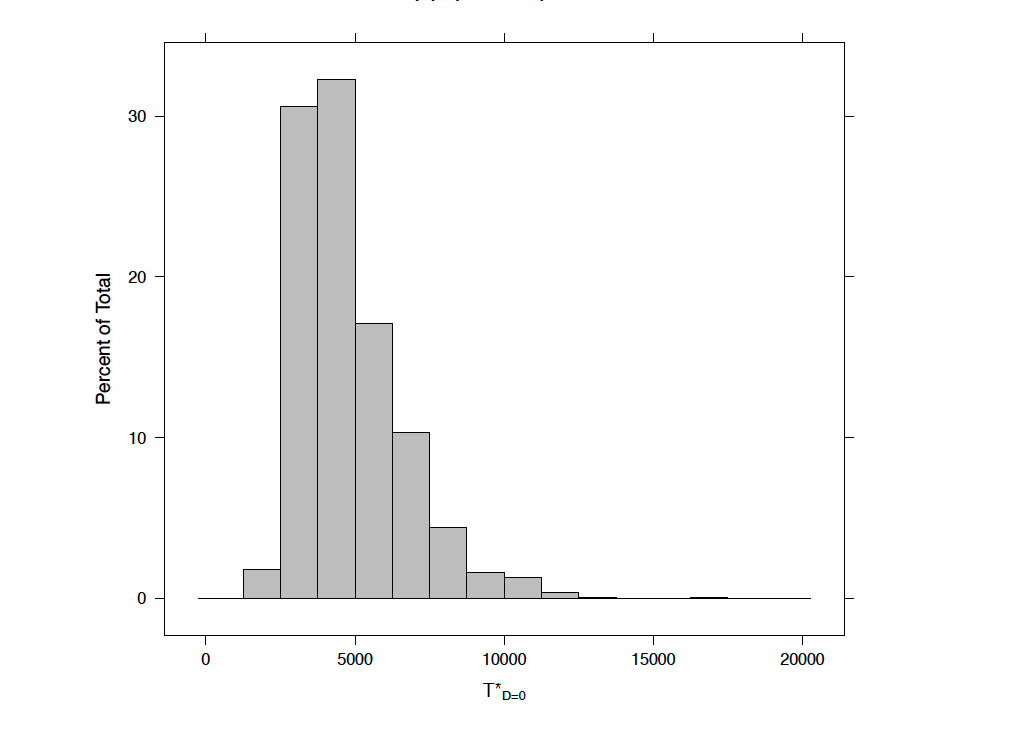}
\caption{\small Distribution of the time $\tau_\Gamma$   for $Y$ to reach $\Gamma$ starting from $y_0=(\frac13,\frac13)$ ($\tau_\Gamma$ is unscaled, $N=1000$).}
\label{fig:46}
\end{figure}

We evaluated the behaviour of the process in its initial phase, by estimating the probability that the process reaches its first extinction before reaching the coexistence balance line $\Gamma=\{d=0\}$ and by estimating the time until it reaches $\Gamma$. For $\tau_\Gamma=\inf\{t:d_t=0\}$ estimates for $P(\tau_\Gamma<\tau^e)$ were calculated for eight different starting points in the right half of the triangle $S$: $(d_0,m_0)$=(0.334,0.334), (0.48,0.48), (0.495,0.495), (0.5,0.5), (0.5,0.020), (0.5,0.01), (0.97,0.01), (0.985,0.005). Using $1000$ runs the estimates showed $P(\tau_\Gamma<\tau^e)< 0.01$ for each starting point except for: $y_0=(0.495,0.495)$ when this probability is $\approx 0.05$ and $y_0=(0.985,0.005)$ when it is $\approx 0.025$. As long as the staring point is at distance $>0.05$ from one of the absorbing corner points of $S$ the pull towards the coexistence line prevails. The empirical probability distribution of the time $\tau_\Gamma$ using  $(d_0,m_0)$=(0.334,0.334)  is shown in Figure~\ref{fig:46}. Note that its upper bound is significantly smaller than $N^2=10^6$ indicating that the initial  phase is negligible relative to the rest of the stochastic dynamics of the process.

\begin{figure}
\centering
\includegraphics[totalheight=0.225\textheight,width=0.725\textwidth]{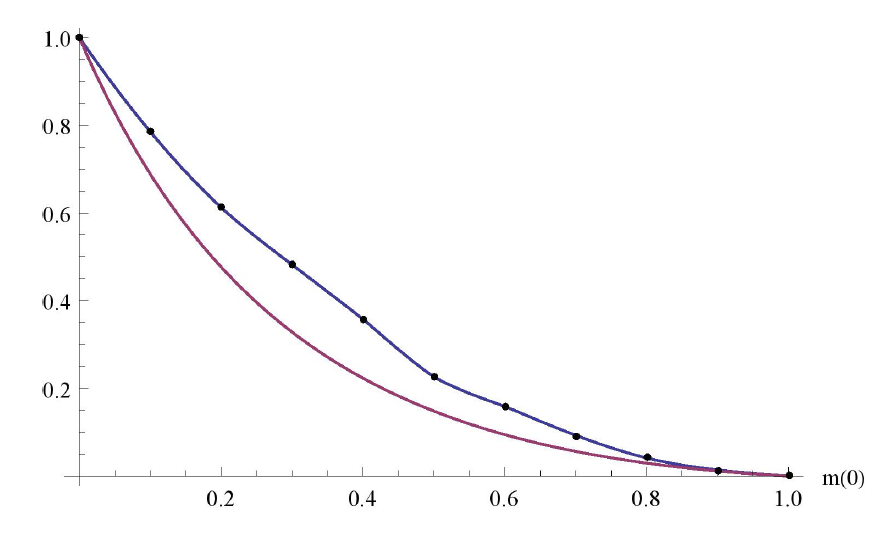}
\caption{\small Graph of $p_M$ as a function of the starting point $x=(0,m_0)$ of the process (with $m_0$ on the $x$-axis). Purple line is based on values calculated from explicit formula (for $1000$ different values of $m_0$); blue line connects values (for $10$ different values of $m_0$) estimated from simulations (using $1000$ runs of the original unscaled process $Y$).}
\label{fig:pm}
\end{figure}

\begin{figure}
\centering
\includegraphics[totalheight=0.225\textheight,width=.85\textwidth]{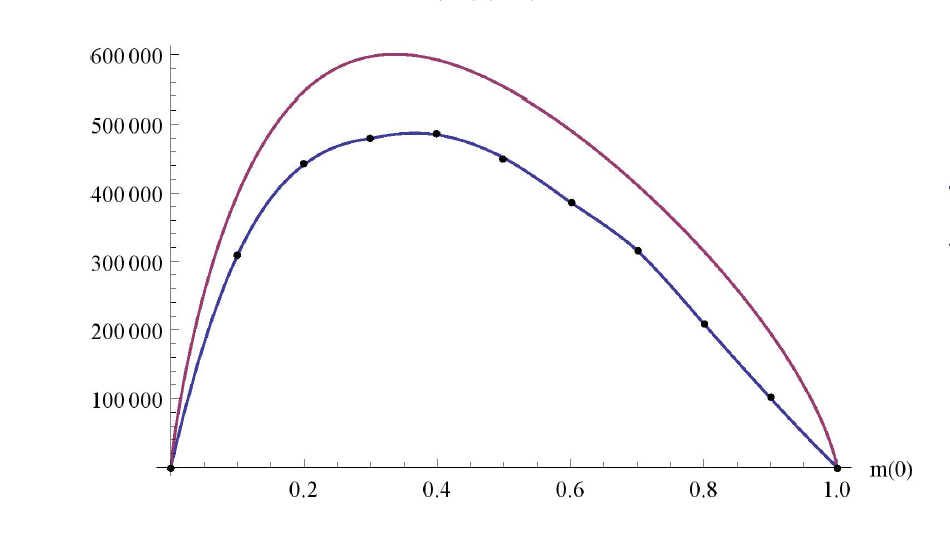}
\caption{\small Graph of $E_x(\tau)$ as a function of the starting point $x=(0,m_0)$ of the process ($m_0$ on the $x$-axis, $N^2E_x(\tau)$ on the $y$-axis for $N=10^3$). Purple line is based on explicit formula (for $1000$ different values of $m_0$); blue line connects values (for $10$ different values of $m_0$) from  simulations ($1000$ runs of unscaled $Y$).}
\label{fig:Etau}
\end{figure}

The dependence of the probability $p^e_M$ and the mean time to the first extinction $E_{y_0}(\tau^e)$ on the value of the starting point for the process $y_0=(d_0,m_0)$ is shown in Figures~\ref{fig:pm} and~\ref{fig:Etau}. The simulations were done using population size $N=1000$ and ten equally spaced values for $y_0=(0,0), (0,0.1), (0,0.2), \dots, (0.09)$, with estimates based in $1000$ runs. Probability of first extinction at $M$, as intuitively expected, decreases with increase in the initial proportion of $M$; in a nonlinear way. The time to extinction at first increases then decreases with increase in the initial proportion of $M$, in a concave way with a maximum at approximately $m_0=0.334$. 

\subsection{Model in a Fluctuating Environment}

In our model the random environment was assumed to be stationary in time, so that the probability distribution of state {\bf c} versus {\bf h} stays constant as the population evolves. In other words the model is embedded in an environment whose values at each replacement event are i.i.d.  with the given distribution $(q,1-q)$. In addition to a time homogeneous random environment this also approximates well the case of a fast fluctuating random environment that equilibriates on the time scale $N$ of competitive dynamics. Our diffusion approximation results then remain valid with $(q,1-q)$ representing the equilibrium distribution of the {\bf c, h} environment.  

In case of slowly fluctuating environment we can approximate the model dynamics by a process of switching diffusions. In case the change in the environment distribution is slower than the time scale $N^2$ of coexistence diffusion in between its changes the competitive dynamics occurs and the species start to diffuse on the line of coexistence balance. When the next change occurs it changes the dynamics leading to new phase of competitive dynamics which results in a different line of coexistence balance. Analysis of dynamics induced by environment fluctuating on a time scale that is between (and including) $N$ and $N^2$ requires a new mathematical approach.


\section{Discussion}

We present a stochastic evolutionary game model and derive rigorous approximations of its long term behaviour in terms of two diffusion processes. Our model incorporates effects of a random environment which determines the fitness of different species and hence influences the competitive dynamics. Without  randomness of the environment stochastic competitive dynamics would lead only towards fixation of the species that is the fittest in that environment. Deterministic competitive dynamics in a random environment would also result only in fixation of the species which is more favoured by the environment (only in the symmetric case would it end in coexistence). In the stochastic  random environment model  the fixation outcome is not predetermined, and moreover, the dynamics goes through a long transient state of coexistence. 

We make predictions for the outcomes of the model, determining probabilities of extinction of species and the expected time to the first extinction event. In order to do so we derive rigorous stochastic approximations characterizing the long term behaviour of the model. Initially its dynamics closely follows  approximately deterministic competition between the species (its mean dynamics). This phase is relatively short and results in a region of 'coexistence balance' between all the species, in the neighbourhood of which it remains for the rest of the time. In the next, longer, phase the dynamics fluctuates within the region of coexistence, and can be approximated by a diffusion process (whose mean and variance coefficients are explicitly derived). Depending on which boundary of the coexistence region this diffusion reaches first, there may be a subsequent even longer phase of coexistence of  two remaining species. In this case the fixation outcome is determined when this latter randomly perturbed dynamical system first reaches one of its boundary points. 

Our results provide a rigorous dimension reduction for the model, which captures the full stochastic nature of its behaviour. It allows one to make calculations and predictions for the model which would analytically be intractable in the original model, with errors that go to zero with increasing population size. The approximating process is a one dimensional diffusions on regions of coexistence for which simple integration yields analytical expressions for the probabilities and expected times. We illustrate the proximity of our approximation to the original stochastic model using simulation results. 

With a simple stochastic model we expose a new approach to long term transient coexistence, which does not fit into the mathematical framework that has so far been considered for analysis of stochastic evolutionary games. We emphasize the usefulness of approximating by degenerate diffusions, providing a rigorous model reduction that accurately captures the long term stochastic behaviour. 



\section{Appendix: Proofs of Propositions}

\subsection{Proof of Propositions~\ref{prop:ODE} and~\ref{prop:SDE1}}
We will combine the results of Propositions~\ref{prop:ODE} and~\ref{prop:SDE1} as they can be shown jointly using the rescaled process $y^{N^2}$ started from a convergent sequence of initial points in $U=S\setminus \{m=0\}$. The proof follows from a very general result of \cite{Katz91} for semimartingales which can be written in terms of a large drift term with a stable manifold of fixed points and a well-behaved fluctuation term. This result applies to the long time behaviour of random perturbations of a dynamical system with a small noise term, provided that the dynamical system has an asymptotically stable manifold of fixed points and is started from a point in its domain of attraction. 

A direct proof of Proposition~\ref{prop:ODE} can be shown using a strong approximation theorem for density dependent population processes established in \cite{Kurtz78} (see also \cite{Kurtz87}). Theorem 2.2 in \cite{Kurtz78} proves that, under appropriate conditions, density dependent Markov chains such as $y^N$ converge to a deterministic process $\bar y$ which is a solution  of the ordinary differential equation \eqref{eq:2dode}. 
The conditions of this theorem are met because the jumps of the Markov process $y^N$ are bounded in size by 2 and because the infinitesimal change in the mean is Lipschitz continuous on any compact set $K_\epsilon \subset U$. 

Once we have the result of Proposition~\ref{prop:ODE}, it is easy to show that then the initial value for the next proposition is met, i.e. for any $t<\tau_\epsilon$ we have $P(\|y^{N^2}_{t/N}-\Gamma\|\le N^{-\delta})\to 1$. A direct proof of Proposition~\ref{prop:SDE1} can be made using the fact that the function $\varphi(y)=(d-(2q-1))^2$ is a stochastic Lyapunov function for the process over the space $S\setminus\{(-1,0),(1,0)\}$, namely there exist constants $0<c_1,c_2<\infty$ such that $\forall y\in K$ in compact $K\subset U$
\begin{eqnarray*}\begin{aligned}\partial_tE\big(\varphi(y^{N^2}_t)|y^{N^2}_t=y\big)&\le -c_1\varphi(y)^2+c_2,\\ c_1&=\inf_{y\in K}\frac{1-m-d^2}{1-d^2},\, c_2=\sup_{y\in K}a_{dd}(y)\end{aligned}\end{eqnarray*}
Further estimates 
together with Doob's maximal inequality for martingales can be used to show that $\varphi(y^{N^2}_{t\land \tau_\epsilon})\to 0$ as $N\to\infty$.

In addition to the separate arguments above we show how both Propositions follow from Katzenberger's powerful results. The rescaled process $(y^N_t)_{t\ge 0}$ has infinitesimal change in the mean $b=(b_d,b_m)=E[\Delta y^N_t|y^N_t]$ and in the (co)variance $a=[(a_{dd},a_{dm})^t,(a_{md},a_{mm})^t]=NVar[\Delta y^N_t|y^N_t]$ given by \eqref{eq:moments}. This implies infinitesimal changes in mean $E[\Delta y^{N^2}_t|y^{^2}N_t]=Nb$ and in the variance $Var[\Delta y^{N^2}_t|y^{N^2}_t]=a$ for the longer term process $y^{N^2}$ (the factor of $N$ is what makes this a process with a large drift term). From expressions in \eqref{eq:moments} it is clear that $b$ is differentiable on $U$ with a bounded derivative on any compact set $K_\epsilon \subset U$. The set of fixed points $b^{-1}(0)=\Gamma$ is one dimensional and the whole set $S\setminus\{(-1,0),(1,0)\}$ is a domain of attraction for $\Gamma$.
The matrix gradient of the vector field $b$ for any $y\in\Gamma$ 
\begin{eqnarray*}
\begin{aligned}
\nabla b(2q-1,m)=\left[\Big(-\frac{1-m-(2q-1)^2}{1-(2q-1)^2}, 0\Big)^t, \Big(\frac{(2q-1)m}{1-(2q-1)^2},0\Big)^t\right]
\end{aligned}
\end{eqnarray*}
has the non-zero eigenvalue $\lambda=-\frac{1-m-d^2}{1-d^2}<0$. Since our state space for is bounded and its jumps are bounded all the other conditions needed to apply Katzenberger's result are satisfied. Theorem 6.2 of \cite{Katz91} implies only the result of the Proposition~\ref{prop:SDE1}, however Theorem 6.3  provides a joint result for both the behaviour in the first phase as in Proposition~\ref{prop:ODE} as well as the behaviour in the longer subsequent phase as in Proposition~\ref{prop:SDE1}. Furthermore, within this framework it is also shown in Theorem 6.1 that the limit point of the stopping times $\tau^N_\epsilon$ from Proposition \ref{prop:ODE} is as $N\to\infty$ almost surely  greater than or equal to the stopping time $\tau$. This implies that the first extinction or fixation event lies on the line $\Gamma$ rather than in either of the excluded corners $\{(-1,0),(1,0)\}$.

%

\begin{acknowledgements} This research was supported by NSERC (Natural
Sciences and Engineering Research Council of Canada) Discovery Grant \# 06573-2015.\\
\end{acknowledgements}

\

\end{document}